\documentclass[11pt]{amsart}
\usepackage{amsmath,amsfonts,amssymb,amsthm,epsfig,color}
\newcommand{\bb}{\mathbb}

\newcommand{\z}{\mathbf z}
\newcommand{\C}{\bb C}
\newcommand{\h}{\bb H}
\newcommand{\Z}{\bb Z}
\newcommand{\R}{\bb R}
\newcommand{\N}{\bb N}

\newcommand{\Q}{\bb Q}

\newcommand{\G}{\bb G}
\newcommand{\s}{\bb S}

\newcommand{\f}{\mathfrak}

\newcommand{\seq}{_{n=0}^{\infty}}
\newtheorem{Theorem}{Theorem}
\numberwithin{Theorem}{section}
\newtheorem{Cor}[Theorem]{Corollary}
\newtheorem{Prop}[Theorem]{Proposition}
\newtheorem{lemma}[Theorem]{Lemma}

\newtheorem*{lemma*}{Lemma}

\newtheorem*{theorem*}{Theorem}

\numberwithin{equation}{section}
\begin{document}
\title{Logarithm laws for unipotent flows, II}
\author{J.~S.~Athreya and G.~A.~Margulis}
\email{jathreya@gmail.com}
\email{margulis@math.yale.edu}
\address{Department of Mathematics, University of Illinois, Urbana, IL}
\address{Department of Mathematics, Yale University, New Haven, CT}
\thanks{J.S.A. supported by NSF grant DMS
    0603636.}
  \thanks{G.M. supported by NSF grant DMS 0801195} 
  \date{4/8/2010}
\begin{abstract}
We prove analogs of the logarithm laws of Sullivan and Kleinbock-Margulis in the context of unipotent flows. In particular, we prove results for horospherical actions on homogeneous spaces $G/\Gamma$. We describe some relations with multi-dimensional diophantine approximation.
\end{abstract}

\maketitle


\section{Introduction}

Two important dynamical systems on non-compact manifolds are the geodesic and horocycle flows on the unit tangent bundle of a finite-volume non-compact hyperbolic surface. Both of these flows are known to be ergodic, and thus, generic orbits are dense. A natural question is to understand the behavior of excursions of trajectories into the cusp(s). 

For geodesic flows, the statistical properties of these excursions were first studied in~\cite{Sullivan} by Sullivan (in the context of finite volume hyperbolic manifolds) and later, in the more general context of the actions of one-parameter diagonalizable subgroups on non-compact finite-volume homogeneous spaces, by Kleinbock-Margulis. In~\cite{KM}, they proved the following result:

\begin{Theorem}\label{loglawzero}(\cite{KM}, Theorem 1.7 and Prop 5.1)  Let $G$ be a connected semisimple Lie group without compact factors, $\f{g}$ its Lie algebra, $\Gamma \subset G$ an irreducible non-uniform lattice, $K$ a maximal compact subgroup, and $d(\cdot,\cdot)$ a distance function on $G/\Gamma$ determined by a right-invariant Riemannian metric on $G$ bi-invariant under $K$. Let $\mu$ denote the measure on $G/\Gamma$ determined by Haar measure on $G$. Let $\f{a} \subset \f{g}$ be  a Cartan subalgebra, $0 \neq \mathbf \mathbf z \in \f{a}$, and $a_t = \exp(t\mathbf z)$. Then there exists a $k = k(G/\Gamma, d)>0$ such that $\forall y$,

\begin{itemize}

\item  $\exists C_1, C_2 >0$ such that for all $t >0$, \begin{equation}\label{cusp} C_1 e^{-kt} \le \mu(x \in G/\Gamma: d(x, y) > t) \le C_2 e^{-kt}.\end{equation}

\item  For $\mu$-a.e. $x$, \begin{equation}\label{gtloglaw}\limsup_{t \rightarrow \infty} \frac{d(a_t x, y)}{\log t} = 1/k.\end{equation}

\end{itemize}

\end{Theorem}

\vspace{.1in}

In this paper, we prove results similar to equation (\ref{gtloglaw}) for several classes of \emph{unipotent} actions. This paper is a sequel to~\cite{AMarg1}, where we considered the case of unipotent flows on the space of lattices. Subsequently, there has been significant activity in the setting of unipotent logarithm laws, for example the papers~\cite{Athreya3, AthreyaPaulin, KelMo}. 

Our results can broadly be divided into two categories:

\begin{enumerate}

\item \emph{Horospherical actions}. We prove a result (Theorem~\ref{horoactions}) on the excursions of orbits of large subsets of horospherical subgroups. We obtain lower bounds for specific orbits. 

\item \emph{Almost everywhere results for flows}. This result (Theorem~\ref{prob}) applies in the most general situation of one-parameter unipotent flows on symmetric spaces, and uses probabilistic techniques (generalized Borel-Cantelli lemmas) and exponential decay of matrix coefficients.

\end{enumerate}

\subsection{Organization:} The paper is organized as follows: in \S\ref{results}, we state our main results. In \S\ref{tori}, we collect technical results on tori and divergent trajectories required for our proofs. In \S\ref{divgeod}, we use these technical results to prove our main theorem on horospherical actions, as well as related corollaries on hyperbolic surfaces. Finally, in \S\ref{borelcantelli}, we prove our probabilistic results.

\subsection{Acknowledgements:} We would like to thank Herbert Abels, Manfred Einsiedler, Alex Eskin, Dmitry Kleinbock, Yair Minsky, Nimish Shah, Yitwah Cheung, Roger Howe, and Barak Weiss for useful discussions.

\section{Statement of results}\label{results}

\subsection{Horospherical actions}\label{hororesults}

Let $G$ be a connected semisimple Lie group without compact factors, and $\Gamma \subset G$ be an irreducible non-uniform lattice. Let $\mu$ denote the probability measure on $G/\Gamma$ arising from Haar measure on $G$. Let $A$ be a maximal connected $\Q$-diagonalizable subgroup, and $\{a_t\}_{t \in \R}$ be a one-parameter subgroup of $A$, and let
\begin{equation}\label{horosphere}H : = \{h \in G: a_{-t} h a_t \rightarrow_{t \rightarrow +\infty} 1\}\end{equation}

\noindent be the \emph{expanding} horospherical subgroup associated to $\{a_t\}$. 

Given $x_0 \in G/\Gamma$, Dani (\cite{Dani}, Theorem 1.6)  proved that if $\{a_{-t} x_0\}_{t \geq 0}$ is non-divergent (i.e., there is a compact set $C \in G/\Gamma$ and a sequence of times $t_n \rightarrow +\infty$ such that $a_{-t_n} x_0 \in C$), that $Hx_0$ is dense in $G/\Gamma$. Our aim is to give a more quantitative version of this result, with regards to visits to neighborhoods of $\infty$. Let $B \subset H$ be a non-empty, bounded, open subset. Set 
\begin{equation}\label{btdef} B_t := a_{\log t} B a_{-\log t}.\end{equation}

\noindent This forms an expanding family of subsets of $H$. 

Let $d_X$ denote a right-invariant metric on $G$ arising from the Riemannian metric induced by the Killing form. If $G^{\prime} \subset G$ is a sugbroup of $G$, we let $d_{G^{\prime}}$ denote the induced distance function on $G^{\prime}$. We let $d_{G/\Gamma}$ denote the induced distance function on $G/\Gamma$. We will drop the subscripts when it is clear on which space we are measuring distances. We will study the behavior of the excursions of $B_t x_0$ away from compact sets by investigating the asymptotic behavior of the quantities 
\begin{equation}\label{eq:beta:def}\beta_t(x_0) : = \sup_{b \in B_t} d_{G/\Gamma}(bx_0, x_0).\end{equation}

\noindent Since $H x_0$ is dense  for all $x_0$ such that $\{a_{-t} x_0\}_{t\geq0}$ is non-divergent, we have $$\limsup_{t \rightarrow \infty} \beta_t(x_0) = \infty$$\noindent for such $x_0$. Our main result is about the rate of these excursions.

To formulate our results, we need a little more notation. If the $\R$-rank of $G$ is at least $2$, we can assume, by the Arithmeticity Theorem (\cite{Margulis}, Chapter IX), that $G = \G(\R)^{\circ}$ and $\Gamma = \G(\Z)$, where $\G$ is a semisimple algebraic $\Q$-group. Let $\s$ be a maximal $\Q$-split torus. Without loss of generality, we can assume $A = \s(\R)^{\circ}$. Let $\|\cdot\|$ denote the norm on $A$ induced by the Killing form. We can write $a_t = \exp{t \mathbf z}$, with $\mathbf z \in \mathfrak{a}$. if the $\R$-rank of $G$ is equal to $1$, then there is (up to scaling and conjugation) a unique $1$-parameter subgroup, which we again can write as $a_t = \exp{ t \mathbf{z}}$ for $z \in \mathfrak{a}$, where $\mathfrak{a}$ is the Lie algebra to $\{a_t\}$. \medskip
%
%
%
\medskip

Given $x_0 \in G/\Gamma$, let 

$$\omega^{-} : = \omega^-(x_0, a_t, d, \Gamma) := \limsup_{t \rightarrow +\infty} \frac{d_{G/\Gamma}(a_{-t}x_0, x_0)}{t}.$$

\begin{Theorem}\label{horoactions} Fix notation as above. Let $a_t = \exp(t \mathbf z)$, $\mathbf z \in \mathfrak{a}$. Let $\nu = \|\mathbf z\|$.

\begin{equation}\label{horoupper} \limsup_{t \rightarrow \infty} \frac{ \beta_t(x_0)}{\log t} \le \nu + \omega^-\end{equation}

\noindent If $\{a_{-t} x_0\}_{t \geq 0}$ is non-divergent, then

\begin{equation}\label{horolower}\limsup_{t \rightarrow \infty} \frac{ \beta_t(x_0)}{\log t} \geq \nu\end{equation}

\end{Theorem}

We will prove this theorem in \S\ref{divgeod}. Combining this result with Theorem~\ref{loglawzero}, which implies that $\omega^-(x_0) = 0$ for $\mu$- almost every $x_0 \in G$, we have the following corollary:

\begin{Cor}\label{horoactionloglaw} Fix notation as above. $$ \limsup_{t \rightarrow \infty} \frac{ \beta_t(x_0)}{\log t} = \nu$$  for $\mu$-almost every $x_0$.\end{Cor}

\vspace{.1in}
\newpage
\noindent\textbf{Remarks:} 

\begin{itemize}

\item It is initially somewhat surprising that the typical horospherical excursion should have the same order as the typical excursion for $a_t$. What we will show in the proof is that the behavior of the horosphere is governed by divergent, non-typical, $a_t$-trajectories.

\medskip

\item (\ref{horoupper}) follows relatively easily from the triangle inequality, whereas (\ref{horolower}) requires a more detailed analysis of divergent $\{a_t\}$-trajectories.

\medskip

\item We expect that our results will hold for general \emph{norm-like pseudometrics}, as defined in~\cite{Abels-Margulis}.

\end{itemize}

\vspace{.1in}

\subsection{Hyperbolic surfaces}\label{hypsurface}

Specializing to $G = SL(2, \R)$ with $H = \{h_s\}_{s \in \R}$, where \begin{equation}\label{horodef} h_s =\left(\begin{array}{cc} 1 & s \\ 0 & 1
\end{array}\right).\end{equation}\noindent we have that \begin{equation}\label{geodef} a_t =\left(\begin{array}{cc} e^{t/2} & 0 \\ 0 & e^{-t/2}
\end{array}\right).\end{equation}  We take $B = \{h_s\}_{s \in (0, 1)}$, and so $B_t = \{h_s\}_{s \in (0, t)}$ ($a_t h_s a_{-t} = h_{se^t}$), and obtain a sharp result for the horocycle flow on the unit tangent bundle of a general non-compact finite volume hyperbolic surface. Let $\Gamma \subset SL(2,\R)$ be a non-uniform lattice. Let  $d$ denote distance on the hyperbolic surface $S = \h^2/\Gamma$ ($\h^2$ denotes the upper-half plane with constant curvature $-1$), and $p: M \rightarrow S$ be the natural projection from $M =SL(2, \R) /\Gamma$.

\begin{Cor}\label{horosl2} Let $H = \{h_s\}_{s \in \R}$.  Fix $y \in S$.Then for all $x \in S$, almost all $\tilde{x} \in p^{-1}(x)$, 
\begin{equation}\label{hyphorolog}\limsup_{s \rightarrow \infty} \frac{d(p(h_s \tilde{x}), y)}{\log s} = 1.\end{equation}
Moreover, for all $\tilde{x} \in M$ such that $H\tilde{x}$ is not closed, 
\begin{equation}\label{horogeq1}\limsup_{s \rightarrow \infty} \frac{d(p(h_s \tilde{x}), y)}{\log s} \geq 1.\end{equation}

\end{Cor}
\bigskip
\noindent The following proposition shows that while (\ref{hyphorolog}) holds for almost every point, the inequality in (\ref{horogeq1}) is strict for a (topologically) large set of points:

\begin{Prop}\label{ubprop} Let $\Gamma \subset SL(2, \R)$ be a non-uniform lattice, $H = \{h_s\}_{s \in \R}$ as in equation~\ref{horodef}. Let $y \in \h^2/\Gamma$. Let $$B = \{x \in SL(2, \R)/\Gamma: \limsup_{t \rightarrow \infty} \frac{d(p(h_s x), y)}{\log s} = 2\}.$$\noindent $B$ contains a dense set of second Baire category.\end{Prop}
\bigskip
\noindent\textbf{Remarks:}\begin{itemize}

\item Note that in the metric on $\h^2$, $d(p(h_s), i) = 2 \log s$ (where by abuse of notation, $p: SL(2, \R) \rightarrow \h^2 = SO(2)\backslash SL(2, \R)$ is the projection $p(g) = SO(2)g$), so $2$ is the maximum value this $\limsup$ can attain. In fact, for any sequence $r_n \rightarrow \infty$ in $SL(2, \R)$, $d(p(r_n), i) \approx 2 \log |r_n|$ ($\approx$ means the ratio goes to $1$), where $|g|$ is the supremum of the matrix entries of $g$.

\medskip

\item By (\ref{horoupper}), the set $B$ must consist of trajectories which diverge at rate $1$ under $a_{-t}$, that is, they must satisfy $$\limsup_{t \rightarrow \infty} \frac{d(p(a_{-t} \tilde{x}), y)}{t} = 1.$$

\medskip

\item In~\cite{AMarg1}, we consider the special case of $\Gamma = SL(2, \Z)$, and obtain several connections to Diophantine approximation. Further results can be found in~\cite{Athreya3}, in which precise conditions for this lim sup to take on certain values for general $SL(2, \R)/\Gamma$ are given, and in~\cite{KelMo} where results are obtained for quotients of products of $SL(2,\R)$ and $SL(2, \C)$.

\end{itemize}

\subsection{Upper and lower bounds}\label{probresults}

We now return to the case of general semisimple Lie groups $G$ and non-uniform lattices $\Gamma$. Now we study the action of \emph{one-parameter} unipotent subgroups on $G/\Gamma$. We have

\begin{Theorem}\label{prob} Fix notation as in Theorem~\ref{loglawzero}. Let $\{u_t\}_{t \in \R} \subset G$ denote a one-parameter unipotent subgroup. Then there is a $0 < \alpha \le 1$ such that for $\forall y$,  $\mu$-a.e. $x$, $$ \limsup_{t \rightarrow \infty} \frac{d(u_t x, y)}{\log t} = \alpha/k $$\end{Theorem}

\medskip

\noindent We prove this theorem in \S\ref{borelcantelli}. \medskip

\noindent\textbf{Remark:} Note that Theorem~\ref{loglawzero} says that if we replace our unipotent subgroup $\{u_t\}$ with a diagonalizable subgroup $\{a_t\}$, we can always take $\alpha =1$. Like Theorem~\ref{loglawzero}, Theorem~\ref{prob} is proved using information on decay of matrix coefficients of the regular representation of $G$ on $G/\Gamma$, and an appropriately adapted version of the Borel-Cantelli lemma. However, the slower decay of matrix coefficients for unipotent flows as compared to diagonalizable flows does not allow us to conclude that $\alpha =1$. It would be very interesting to find examples of unipotent subgroups where $\alpha \neq 1$, though we suspect that such subgroups do not exist.

\medskip

\section{Divergent trajectories}\label{tori}

Recall the notation of \S\ref{hororesults}: $G$ is a connected semisimple Lie group without compact factors, $\Gamma$ an irreducible non-uniform lattice, and $d$ the Riemannian metric arising from the Killing form on $G$. We also use $d$ to denote the metric on $G/\Gamma$. $A$ is a maximal $\Q$-diagonalizable subgroup of $G$, and $\{a_t\}$ a one-parameter subgroup of $A$. Write $a_t = \exp{t \mathbf z}, \mathbf y \in \mathfrak{a}$, and let $\nu = \|y\|$.
We prove the following result concerning $\{a_t\}$ trajectories in $G/\Gamma$.

\begin{Prop}\label{divtraject} Fix $x_0 \in G/\Gamma$.  For all $x \in G/\Gamma$
\begin{equation}\label{nulimsup} \limsup_{t \rightarrow \infty} \frac{d(a_t x, x_0)}{t} \le \nu. \end{equation}

\noindent Moreover, for all $x = g\Gamma$ with $g \in \G(\Q)$, 
\begin{equation}\label{nulimit} \lim_{t \rightarrow \infty} \frac{d(a_t x, x_0)}{t} = \nu. \end{equation}
\end{Prop}

\medskip

\subsection{Reduction Theory} We recall some results from reduction theory. Assume the $\R$-rank of $G$ is greater than $1$. We can assume, as in \S\ref{hororesults}, $G = \G(\R)^{\circ}$ and $\Gamma = \G(\Z)$, where $\G$ is a semisimple algebraic $\Q$-group. Let $\s$ be a maximal $\Q$-split torus in $\G$, and set $A = \s(\R)^{\circ}$. 

Let $\Phi$ be a system of $\Q$-roots associated to $\mathfrak{A}$ and let $\Phi^+$ and $\Phi^s$ be the sets of positive and simple roots respectively. We define the positive \emph{Weyl chamber} 
$$\mathfrak{a}^+ = \{ \mathbf \mathbf z \in \mathfrak{d}: \alpha(\mathbf z) \geq 0 \mbox{ for all } \alpha \in \Phi^s\}.$$

Using the exponential map, we identify it with $A^+ =\exp(\mathfrak{a}^+) \subset A$. Conjugating if necessary, we can assume that $\mathbf z \in \mathfrak{a}^+$, that is, $a_t \in A^+$ for $t >0$. We have the Iwasaswa decomposition $G = KAMU$ (here, $K$ is a maximal compact subgroup, $U$ is unipotent, and $M$ is reductive, with $A$ centralizing $M$ and normalizing $U$). Let $Q \subset MU$ be relatively compact, and for $\tau >0$, define
$$\mathfrak{a}_{\tau} = \{\z \in \mathfrak{a}: \alpha(\mathbf z) \geq \tau \mbox{ for all } \alpha \in \Phi^s\}$$

\noindent We can define a \emph{generalized Siegel set} 
$$S_{Q, \tau} : = K \exp(\f{a}_{\tau}) Q.$$

\noindent For appropriate choices of $Q$ and $\tau$, a finite union of translates of $S_{Q, \tau}$ form a weak fundamental domain for the $\Gamma$-action on $G$. Precisely, we have

\begin{Theorem}[\cite{Leuzinger}, Proposition 2]\label{fundset} Fix notation as above. There is are $Q, \tau$ and $\{q_1, \ldots, q_m\} \in \G(\Q)$ so that $\Omega : = \bigcup_{i=1}^m S_{Q, \tau} q_i$ satisfies

\begin{enumerate}

\item $G = \Omega \Gamma$

\item For all $q \in \G(\Q)$, $\{\gamma \in \Gamma: \Omega q \cap \Omega \gamma \neq \emptyset\}$ is finite.

\end{enumerate}

\end{Theorem}

\medskip

\noindent The finite set ${q_1, \ldots, g_m}$ form a set of representatives for the double coset space $\mathbb{P}(\Q)\backslash \G(\Q) /\Gamma$, where $\mathbb{P}(\Q)$ is a minimal $\Q$-parabolic subgroup. Note that $\mathbb{P}(\R)^{\circ} = AMU$. Fix $x_0$ to be the identity coset in $G/\Gamma$. Let $x = g\Gamma$. Note that, letting $e$ denote the identity in $G$, we have, for any $h \in G$,

$$d(h x , x_0) = \inf_ {\gamma \in \Gamma} d(h g, \gamma)$$

Using Theorem~\ref{fundset}, Leuzinger~\cite[Theorem 1]{Leuzinger} proved that there is a $b \in \overline{A^+}$ (here, $\overline{A^+}$ denotes the closure of the Weyl chamber $A^+$) such for any $\bf y \in \f{a}^+$, with $\| \bf y \| = 1$, $a_t : = \exp(t \bf y)$, any $p \in MU$, any $q_i, 1\le i \le m$, and $\gamma \in \Gamma$, we have
\begin{equation}\label{eq:divergent:rate}
d(a_t b p q_i, p b q_i \gamma) \geq t
\end{equation}

\subsection{Proof of Proposition~\ref{divtraject}} First note that the upper bound (\ref{nulimsup}) follows from the definition of distance on the quotient and the fact that $\|\bf z\| = \nu$. To show the limit result (\ref{nulimit}), we first consider the case when $\R$-rank is at least $2$, so we can apply Theorem~\ref{fundset} and equation~(\ref{eq:divergent:rate}). Since we can write each element $g \in \G(\Q)$ as $g = p q_i \gamma_0$ for $p \in \mathbb{P}(\Q)$ and $\gamma \in \Gamma$, and denoting the bounded error given by the element $b$ by $C$, we have, for all $\gamma \in \Gamma$,
\begin{equation}
d(a_t g, \gamma) \geq \nu t  - C
\end{equation}

\noindent (\ref{nulimit}) follows immidiately.

\subsubsection{$\R$-rank 1} Finally, suppose the $\R$-rank of $G$ is $1$. Applying standard reduction theory~\cite{GR} and the density of orbits of parabolic subgroups (\cite{Mostow}, Lemma 8.5) there is a dense set of points diverging under $a_t$ at rate $\nu = \|z\|$. See also~\cite{Dani, Weiss} for more details on divergent trajectories.

\qed\medskip

\section{Horospherical actions}\label{divgeod}

We fix notation as in \S\ref{hororesults} and \S\ref{tori}.  The proof of Theorem~\ref{horoactions} splits naturally into an upper and lower bound:

\subsection{Lower bound}

\begin{lemma}\label{lower bound} For all $x \in G/\Gamma$ with $\{a_{-t} x\}_{t >0}$ non-divergent, \begin{equation}\label{lbhoro}\limsup_{t \rightarrow \infty} \frac{ \beta_t(x)}{\log t} \geq \nu.\end{equation}\end{lemma}

\noindent\textbf{Proof:} The idea is as follows: given the piece of orbit $B_{e^T} x$, we want to show that it has moved depth $T$ into the cusp. We can write $B_{e^T} x = a_t B a_{-t} x$. If $a_{-t}x$ is non-divergent, we can take some $T$ so that $a_{-t}x$ is in a compact set. Using the fact the forward divergent $\{a_t\}$ trajectories are dense, we can find a divergent trajectory (moving at rate $\nu$) in a `thickening' of the orbit $B a_{-t} x$ in the directions transverse to $H$. Since $a_t$ does not expand the directions transverse to $H$, the divergent trajectory (which will be approximately depth $\nu T$ into the cusp after applying $a_t$) will be near $B_{e^T}x$, so there is some $h \in B_{e^T}$ with $hx$ almost depth $T$ into the cusp, as desired. To make this argument precise, we need to use the following 

\begin{lemma}\label{divtraj} Let $C \subset G/\Gamma$ be compact with non-empty interior, and $\epsilon, \phi >0$. Then there is a $T_{C, \epsilon, \phi}$ such that 
$$\left\{x: d(a_t x, x_0) > (\nu - \phi) t - T_{c, \epsilon, \phi} \mbox{ for all } t>0\right\}$$ 
\noindent is $\epsilon$-dense in $C$.\end{lemma}

\noindent\textbf{Proof:} Note that by Prop~\ref{divtraject},
$$\left\{ x \in G/\Gamma: \exists T(x) \mbox{ such that } d(a_t x, x) >  (\nu - \phi)t - T(x) \mbox{ for all } t>0\right\}$$ 
\noindent is dense in $G/\Gamma$. 

Now let $\epsilon >0$, $C \subset G/\Gamma$ compact. Let $\{ B(x, \epsilon)\}_{x \in C}$ be the cover of $C$ by open $\epsilon$-metric balls. Since $C$ is compact, we can take a finite subcover $\{D_1, D_2, \ldots, D_n\}$, where each $D_i = B(x_i, \epsilon)$. For $1 \le i \le n$, there is a $x_i \in D_i$, $T(x_i) >0$, such that 
$$d(a_t x_i, y) > (\nu - \phi) t - T(x_i).$$\noindent 
Let $T_{C, \epsilon, \phi} = \max_{1 \le i \le n} T(x_i)$. Now, for all $x \in C$, there is an $x_i$ such that $d(x_i, x) < \epsilon$, and for all $1 \le i \le n$, $d(a_t x_i, y) > (\nu - \phi)t- T_{C, \epsilon, \phi}$, so we have our result. \qed\medskip

\vspace{.1in}

Let $H^{-0}$ be the subgroup associated to the neutral/stable directions for $a_t$ ($t>0$). Let $x \in G/\Gamma$ be such that $a_{-t} x$ is non-divergent. Thus, there is a compact $C^{\prime\prime} \subset G/\Gamma$ be compact with a non-empty interior and $t_n \rightarrow \infty$ so that $a_{-t_n} x \in C^{\prime\prime}$ for all $n$. 

We fix one more piece of notation: letting $G^{\prime}$ be a subgroup of $G$, $g_0 \in G^{\prime}$, $r>0$, we let 
$$B_{G^{\prime}}(g_0, r) := \{ g \in G^{\prime}: d_{G^{\prime}}( g_0, g) < r\}.$$\noindent 
Let $\epsilon_1$ be such that for all $\epsilon < \epsilon_1$, there are $\epsilon^{+}, \epsilon^{-}$, 
$$B_G( \epsilon) = B_{H^{-0}} (\epsilon^-)B_H (\epsilon^+).$$

\noindent Let $C^{\prime} = \overline{BC^{\prime\prime}}$. Let $b_0 \in B$, $\epsilon_0 >0$ such that $B_{H} (b_0, \epsilon_0) = B_{H}(\epsilon_0) b_0 \subset B$. There is an $0 < \epsilon < \epsilon_1$, and an $\epsilon_{\prime}$ so that (perhaps shrinking $\epsilon_0$) we can write 
$$B_G(\epsilon) = B_{H^{-0}}(\epsilon^{\prime})B_{H}(\epsilon_0).$$
\noindent Let $C = \overline{B_G(\epsilon) C^{\prime}}$. Now $b_0 x_n \in C^{\prime}$, so $B_G(\epsilon)b_0 x_n \in C$. Shrinking $\epsilon$ if necessary, we have 
$$B_G(\epsilon) b_0 x_n = B(b_0x_n, \epsilon).$$
\noindent Fix $\phi >0$, and let $T = T_{C, \epsilon, \phi}$. There is an $x_n^{\prime} \in B(b_0 x_n, \epsilon)$ so that 
$$d(a_{t_n} x_n^{\prime}, x_n^{\prime}) > (\nu - \phi)t_n - T.$$
\noindent We can write $x_n^{\prime} = h^{-} b_0 x_n$ for $h^{-} \in B_{H^{-0}} (\epsilon^{\prime})$. Now we have 
$$a_{t_n} x_n^{\prime} = h_n^- b_n x,$$
\noindent where $$h_n^- = a_{t_n} h^- a_{-t_n} \in B_{H^{-0}} (\epsilon^{\prime})$$\noindent and 
$$b_n = a_{t_n} b_0 a_{-t_n} \in B_{e^{t_n}}.$$
Thus, we have 
$$d(b_n x, y) \geq d(a_{t_n} x_n^{\prime}, x_n^{\prime}) -\epsilon \geq (\nu -\phi)t_n - T - \epsilon$$
\noindent So 
$$\lim_{n \rightarrow \infty} \frac{d(b_nx, x)}{t_n}  \geq (\nu - \phi)$$
\noindent as $n \rightarrow \infty$,(note that since $x_n^{\prime}$ varies in a compact set, it does not matter in the limit whether we measure distance from $x$ or $x_n^{\prime}$). Thus, 
$$\lim_{n \rightarrow \infty} \frac{\beta_{e^{t_n}} (x)}{t_n} \geq \nu -\phi$$
\noindent which, since $\phi>0$ was arbitrary yields our result. \qed\medskip

\subsection{Upper bound}

\begin{lemma}\label{upper bound} For all $x \in G/\Gamma$, \begin{equation}\label{ubhoro}\limsup_{t \rightarrow \infty} \frac{ \sup_{h \in B_t} d(hx, y)}{\log t} \le \nu + \omega(x).\end{equation}\end{lemma}

\noindent\textbf{Proof:} Let $\epsilon >0$. By the definition of $\omega$, and the boundedness of $B$ for all $t$ sufficiently large, for all $b \in B$, 
$$d(b a_{-\log t} x, x) < (\omega + \epsilon)\log t.$$
\noindent By definition 
$$d(a_{\log t} b a_{-\log t}, b a_{-\log t}) \le \nu \log t.$$  
Combining these two inequalities, and using the triangle inequality, we have, for all $b \in B$ and $t$ sufficiently large, 
$$d(a_{\log t} b a_{-\log t} x, x) < (\omega + \nu + \epsilon) \log t.$$

\noindent Since $\epsilon$ was arbitrary, we have our result.\qed\medskip

\noindent\textbf{Proof of Theorem~\ref{horoactions}:} Combine Lemmas~\ref{upper bound} and ~\ref{lower bound}\qed

\subsection{Hyperbolic Geometry} In this subsection we prove Corollary~\ref{hyphorolog} and Proposition~\ref{ubprop}.

\medskip

\noindent\textbf{Proof of Corollary~\ref{hyphorolog}:} Apply Theorem~\ref{horoactions} to $G = SL(2, \R)$, with $H, B$ and $\{a_t\}$ as in \S\ref{hypsurface}. Note that the Riemannian metric on $\h^2 = K\backslash G$ is coarsely isometric to the normlike metric induced by the norm it induces on $A$. 
\qed\medskip

\noindent\textbf{Proof of Proposition~\ref{ubprop}:}

We need the following lemma, which exploits properties of divergent geodesic trajectories:

\begin{lemma}\label{divhor} Let $E \subset SL(2, \R)/\Gamma$ be open, and $y \in \h^2/\Gamma$. There is a $C = C(E)$ such that for all $T>0$ there is a $z \in E$, $t > T$ such that \begin{equation}\label{ubeps}d(p(h_t z), y)> 2\log t - C\end{equation}\end{lemma}

\noindent\textbf{Proof:} By the density of divergent \emph{geodesic} trajectories there is a $c = c(A)>0$ and a $z \in E$ such that $d(p(g_s z),y) > s - c$ for all $s>0$ (for the rest of this section, we will use the notation $ a_t =\left(\begin{array}{cc} e^{t/2} & 0 \\ 0 & e^{-t/2}
\end{array}\right)$). 

Fix lifts of $y$ and $z$ to a fundamental domain for $\Gamma$ in $\h^2$, call them $y_0$ and $z_0$ ($y_0$ will be a point, $z_0$ will be a point and a unit tangent vector). There will be a horocycle connecting $p(z_0)$ and $p(g_s z_0)$, and as $s \rightarrow \infty$, the inward pointing tangent vector to this will approach the vector $z_0$.

More precisely, suppose without loss of generality $z_0$ is $i$ with the upward pointing tangent vector, i.e., $z_0 = e \in SL(2, \R)$. Then $p(g_s z_0) = p(g_s) = SO(2)g_s$, and if $v_s = r_{\theta_s} \in SO(2)$ is the unit tangent vector (based at $i = p(z_0)$) determining the horocycle connecting $e^{s} i = p(g_s)$ and $i$, we have that $v_s$ approaches the upward pointing tangent vector as $t \rightarrow \infty$, or equivalently $\theta_s \rightarrow 0$. 

In addition, if $t = t_s$ is the time it takes for the horocycle to reach $e^{s} i$, we have $SO(2) g_s = SO(2)h_t r_{\theta_s}$, i.e., there is a $\theta^{\prime}_s$ such that $h_t = r_{\theta^{\prime}_S} g_s r_{\theta_s}$ (this is simply the Cartan (or $KAK$) decomposition). It is an easy calculation that $t_s \approx e^{s/2}$, or equivalently, $s \approx 2 \log t_s$. Thus, for $s >>0$, $r_{\theta_s} \mathbf z \in A$, and 
$$d(p(h_{t_s}r_{\theta_s}z),y) = d(p(g_s z), y) > s -c > 2\log t_s - C,$$
\noindent for some possibly larger $C$.\qed\medskip

To complete the proof of the proposition, define $f_T: SL(2, \R)/\Gamma \rightarrow [0, 2)$ by 
$$f_T(x) = \sup_{2 \le t\le T} \frac{d(p(h_tx), y)}{\log t}.$$
\noindent $f_T(x)$ is increasing in $T$, and bounded, so we can define $f_{\infty}(x) = \lim_{T \rightarrow \infty} f_T(x)$. The $f_T$'s are continuous for $T< \infty$, but $f_{\infty}$ is not. We have 
$$B = \{x: f_{\infty}(x) = 2\} = \bigcap_{k=1}^{\infty} \bigcup_{n=0}^{\infty}\left \{x: f_n(x) > 2-\frac{1}{k}\right\}.$$
\noindent Now, for each $k$, $\bigcup_{n=0}^{\infty}\left \{x: f_n(x) > 2-\frac{1}{k}\right\}$ is dense by Lemma~\ref{divhor}, and open by the continuity of $f_n$. Thus $B$ is a countable intersection of open dense sets, as desired.\qed\medskip

\section{Borel-Cantelli lemmas}\label{borelcantelli}

In this section we prove Theorem~\ref{prob}, using a generalization of the Borel-Cantelli lemma. The classical Borel-Cantelli lemma is as follows:

\begin{lemma}\label{bc}(Borel-Cantelli) Let $\{X_n\}_{n=0}^{\infty}$ be a sequence of $0-1$ random variables, with $P(X_n = 1) =: p_n$. Then, if $\sum_{n=0}^{\infty} p_n < \infty$, 
$$P\left(\sum_{n=0}^{\infty} X_n = \infty\right) = 0.$$ 
If the $X_n$'s are pairwise independent, we have that 
$$P\left(\sum_{n=0}^{\infty} X_n = \infty\right) = 1.$$ if $\sum_{n=0}^{\infty} p_n = \infty.$ 
\end{lemma}

The first statement (non-independent) statement is the `easy half' of this Lemma, and can be derived by simply doing an expectation calculation.

The first example of a logarithm law can be derived from the lemma as follows. Fix $\lambda >0$. Let $\{Y_n\}_{n=0}^{\infty}$'s be independent identically distributed (i.i.d.) exponential random variables with parameter $\lambda$. That is, for any $t >0$,  $$P(Y_n > t) = e^{-\lambda t}.$$ \noindent Let $\{r_n\}\seq$ be a sequence of positive real numbers. Applying Lemma~\ref{bc} to the sequence of random variables  

\begin{equation}X_n : = \left\{
\begin{array}{ll}  1 & Y_n > r_n\\ 0 & \mbox{otherwise},\end{array}\right.\nonumber\end{equation}

\noindent implies $Y_n > r_n$ infinitely often if and only if $\sum_{n=0}^{\infty} e^{-\lambda r_n} = \infty$. As a corollary, one obtains that almost surely 
$$\limsup_{n \rightarrow \infty} \frac{Y_n}{\log n} = 1/\lambda.$$
\noindent To prove Theorem~\ref{prob}, we use the following (relatively standard) generalization of Lemma~\ref{bc} to weakly dependent sequences.

\begin{Prop}\label{genprob} Let $(S, \Omega, P)$ be a probability space (i.e., $\Omega$ is a $\sigma$-algebra of subsets of $S$, and $P: \Omega \rightarrow [0,1]$ is a probability measure). Let $X_n: S \rightarrow \{0,1\}$ be a sequence of $0-1$ random variables on $S$, with $P(X_n = 1) =: p_n$. Also define $p_{i,j} := P(X_iX_j =1)$. Suppose

\begin{enumerate}
\item $\sum_{n=1}^{\infty} p_n = \infty.$
\item There is a function $\psi(m)$ such that  for all $m >0$, $$\sup_{n} |p_{n, n+m} - p_n p_{n+m}| \le \psi(m).$$ 
\item $$\lim_{n \rightarrow \infty} \frac{\sum_{m=1}^{n} \psi(m)(n-m)}{(\sum_{i=1}^n p_i)^2} =0.$$
\end{enumerate}
Then $$P(\sum_{n=0}^{\infty} X_n = \infty) = 1.$$
\end{Prop}

\noindent\textbf{Proof:} Given measurable $X: S \rightarrow \R$, we write $E(X) := \int_{S} X dP$ for the expectation, and $V(X) = E(X^2) - E(X)^2$ for the variance.

Let $J_n = \sum_{i=1}^n X_i$, and $$Y_n = \frac{J_n}{\sum_{i=1}^n p_i} = \frac{J_n}{E (J_n)}.$$\noindent We will show that for any $\epsilon >0$, $P(|Y_n -1| > \epsilon) \rightarrow 0$, which will imply that there is a sequence $n_k$ such that $Y_{n_k} \rightarrow 1$ with probability 1, and thus, that $J_{n_k} \rightarrow \infty$. 

Since $E(Y_n) =1$, it suffices to show that $V(Y_n) \rightarrow 0$, that is, that $Y_n \rightarrow 1$ in $L^2$. Now, 
$$V(Y_n) = \frac{V(J_n)}{E(J_n)^2}.$$
\noindent We have
$$V(J_n)  = V(\sum_{i=1}^n X_i) = \sum_{i=1}^n V(X_i) + \sum_{1 \le i \neq j \le n} Cov(X_i, X_j),$$
\noindent where $\mbox{Cov}(X_i, X_j) : = |p_{i,j} - p_ip_j|$ is the \emph{covariance} of $X_i$ and $X_j$. 

Now, $\mbox{Cov}(X_i, X_j) \le \psi(|j-i|)$ by property (2), so we get that 

\begin{eqnarray}
V(J_n) \le \sum_{i=1}^n p_i + 2\sum_{i=1}^n \sum_{j= i+1}^n \psi(j-i) &=& \sum_{i=1}^n p_i+ 2\sum_{i=1}^n \sum_{m=1}^{n-i} \psi(m) \nonumber\\ &=& \sum_{i=1}^n p_i + 2\sum_{m=1}^n \psi(m)(n-m).
\end{eqnarray}

Dividing by $(\sum_{i=1}^n p_i)^2$, we get that the two right hand terms go to zero (by properties (1) and (3) respectively), and thus, we have our result. \qed \medskip

We would like to apply this result to the context of group actions on homogeneous spaces. Fixing notation as in \S\ref{horoactions}, and given $y \in G/\Gamma$, we define a sequence of functions $Y_n: G/\Gamma \rightarrow \R^+$ by $Y_n(x) = d(u_n x, y)$. Given a sequence of numbers $\{r_n\}_{n \in \N}$, we set 

\begin{equation}X_n(x) : = \left\{
\begin{array}{ll}  1 & Y_n(x) > r_n\\ 0 & \mbox{otherwise.}\end{array}\right.\nonumber\end{equation}

\noindent In order to apply Proposition~\ref{genprob} to our context, we must estimate two quantities:

\begin{enumerate}

\item $\mu(x: d(x, y) > t)$
\item The covariances for the random variables $X_n$. 

\end{enumerate}

The first estimate follows from equation (\ref{cusp}), which yields (since $u_n$ is measure preserving):
$$C_1e^{-kr_n} \le p_n = \mu(x: X_n(x) = 1) \le C_2e^{-kr_n}.$$
In order to estimate the covariances, we must control the matrix coefficients of the sequence $\{u_k\}$ under the regular representation of $G$ on $G/\Gamma$. To do this, we turn once again to~\cite{KM}. The following result is essentially a combination of Proposition 4.2 and Corollary 3.5 from that paper:

\begin{Prop}\label{matrix}\cite{KM} There are constants $C>0$, $0 < \beta <1$ such that for all $n, m \in \N$ $$|p_{n, n+m} - p_ np_{n+m}| \le Cp_n p_{n+m} m^{-\beta},$$ where $p_{i, j} = \mu(x: X_i(x)X_j(x) =1)$.\end{Prop}

\noindent\textbf{Remark:} If we were able to obtain $\beta \geq 2$, we would in fact be able to prove $\alpha =1$ in the statement of Theorem~\ref{prob} following Prop 4.1 in ~\cite{KM}. However, for reasons beyond the scope of this paper, $\beta <2$. 

\medskip

We will not prove Proposition~\ref{matrix} in this paper, instead referring the interested reader to the appropriate sections of~\cite{KM}.

\medskip

\noindent\textbf{Proof of Theorem~\ref{prob}:} Let $r_n > \frac{1}{k} \log n$. Then $p_n$ is summable, so for almost all $x$, $X_n = 1$ only finitely often, yielding our upper bound. For our lower bound we apply Proposition~\ref{genprob} to our sequence $X_n$, with $\psi(m) = m^{-\beta}$. It is a simple calculation that for any $\gamma < \beta/2$, setting $r_n = \frac{\gamma}{k} \log n$ will yield:
$$\lim_{n \rightarrow \infty} \frac{\sum_{m=1}^{n} \psi(m)(n-m)}{(\sum_{i=1}^n p_i)^2} =0.$$
\noindent Using Proposition~\ref{genprob}, we have, for $\mu$-a.e. $x$,  
$$\beta/2k \le \frac{\limsup_{t \rightarrow \infty} d(u_t x, y)}{\log t} \le 1/k.$$
\noindent Finally, note that  
$$\limsup_{t \rightarrow \infty} \frac{d(u_t x, y)}{\log t}$$
\noindent is a measurable $u_t$-invariant function on $G/\Gamma$. Thus, if the $u_t$-action is ergodic, it must be constant almost everywhere. If $u_t$ is not ergodic, it must act trivially in some factor of $G$ by the Moore ergodicity theorem~\cite{Moore}, and thus we can reduce to the ergodic case.\qed

\end{document}